\documentstyle[12pt]{article}
\newcommand{\fr}{ \frac}
\newcommand{\beq}{\begin{equation}}
\newcommand{\eeq}{\end{equation}}
\newcommand{\beqa}{\begin{eqnarray}}
\newcommand{\eqa}{\end{eqnarray}}
\newcommand{\lb}{\label}
\newcommand{\r}{\ref}
\newcommand{\G}{\Gamma_q}
\newcommand{\p}{\psi_q}
\newcommand{\ar}{\rightarrow}

\begin{document}
\begin{flushright}
FGI-99-4 \\
math.QA/9812145
\end{flushright}

\begin{center}
{\Large {\bf Quantum Analogue of the Neumann Function of Integer Order}}
\end{center}

\vspace{1cm}

\begin{flushleft}
H. Ahmedov$^1$  and I. H. Duru$^{2,1}$
\vspace{.5cm}

{\small 
1. Feza G\"ursey Institute,  P.O. Box 6, 81220,  \c{C}engelk\"{o}y, 
Istanbul, Turkey 
\footnote{E--mail : hagi@gursey.gov.tr and duru@gursey.gov.tr}.

2. Trakya University, Mathematics Department, P.O. Box 126, 
Edirne, Turkey.}
\end{flushleft}

\vspace{1cm} \noindent 
{\bf Abstract}: q-Neumann function of integer order $N_n (x;q)$ is obtained
and some of its properties are given. q-Psi function which is used in
deriving $N_n (x;q)$ is also introduced and some of its properties are
presented.

\begin{center}
March 1999
\end{center}

\vspace{1cm} \noindent
{\large {\bf 1. Introduction}}

\vspace{2mm} \noindent
The Hahn-Exton q-Bessel functions $J_\nu (x; q)$ which are closely connected
to the quantum group of the plane motions are well studied \cite{1,2,3,4}.
Note that for $\nu =n$ is integer the functions $J_n(x; q)$ and $J_{-n}(x; q)
$ are not independent of each other. To our knowledge the quantum analogues
of the Bessel functions of integer order which are independent of $J_n(x;q)$
have not been addressed.

The main purpose of this note is to introduce the second independent
solution $N_n (x;q)$ of the q-Bessel difference equation. This function
possess the same recurrence relations as the Hahn-Exton q-Bessel function
and in $q\rightarrow 1^{-}$ becomes Neumann function of the order $n$. We
call it the q-Neumann function of order $n$. This solution is non regular at 
$x=0$. It is well known that non-regular solutions are important, since the
Green functions are given in terms of them. For example the q-Legendre
function of the second kind shows up as the Green function on the quantum
sphere \cite{5}. The recently obtained Green function on the quantum plane
is in fact the superposition of the Hahn-Exton q-Bessel and q-Neumann
functions of the order $0$ \cite{6}.

Since the classical Neumann function of integer order $N_n(x)$ is obtained
by taking the derivatives of the Bessel functions $J_\nu (x)$ and $J_{-\nu}
(x)$ with respect to the order $\nu$, it involves the psi functions.
Therefore to derive the q-Neumann function of integer order we first have to
have q-psi function in hand.

The Section 2 is devoted to the introduction of the q-psi function and some
of its properties which are employed to derive the q-Neumann function of
integer order.

In Section 3 we obtained the q-Neumann function of integer order and
presented some relations involving it.

\vspace{2cm} \noindent
{\large {\bf 2. q-Psi Function }}

\vspace{2mm} \noindent
We define the q-psi function as 
\begin{equation}
\psi _q(\nu )=\frac d{d\nu }\log \Gamma _q(\nu ),
\end{equation}
where the q-gamma function $\Gamma _q(\nu )$ is defined by ( $0<q<1$ ) 
\begin{equation}
\label{g}\Gamma _q(\nu )=(1-q)^{1-\nu }\prod_{l=1}^\infty \frac{1-q^l}{%
1-q^{\nu +l}}.
\end{equation}
Many properties of the gamma function were derived by Askey \cite{a}. It is
obvious from (\ref{g}) that $\Gamma _q(\nu )$ has poles at $\nu =0,\ 1,\ 2,\
.\ .\ .$. The residue at $\nu =-n$ is 
\begin{equation}
\label{rg}\lim _{\nu \rightarrow -n}(\nu +n)\Gamma _q(\nu )=(-1)^n\frac{%
(q-1)q^{-n(n+1)/2}}{\log q\Gamma _q(n+1)}.
\end{equation}
The explicit form of $\psi _q(\nu )$ is 
\begin{equation}
\psi _q(\nu )=-\log (1-q)+\log q\sum_{l=0}^\infty \frac{q^{\nu +l}}{1-q^{\nu
+l}}
\end{equation}
The recurrence relations and asymptotic conditions satisfied by this
function are 
\begin{equation}
\label{r}\psi _q(\nu +n)=\psi _q(\nu )-\log q\sum_{l=0}^{n-1}\frac{q^{\nu +l}%
}{1-q^{\nu +l}},
\end{equation}
\begin{equation}
\psi _q(\nu -n)=\psi _q(\nu )+\log q\sum_{l=0}^n\frac{q^{\nu -l}}{1-q^{\nu
-l}}
\end{equation}
and 
\begin{equation}
\lim _{\nu \rightarrow \infty }\psi _q(\nu )=-\log (1-q);\ \ \ \lim _{\nu
\rightarrow -\infty }\psi _q(\nu )=\infty .
\end{equation}
$\psi _q(\nu )$ has poles at $\nu =0,\ 1,\ 2,\ .\ .\ .$ with the residue 
\begin{eqnarray}\lb{rp}
\lim_{\nu\ar -n} (\nu +n)\p (\nu ) & = & 
\log q\lim_{\nu\ar -n} (\nu +n)
\sum_{l=0}^\infty \fr{q^{\nu + l}}{1-q^{\nu+l}} \nonumber \\ 
& = & \log q \lim_{\nu\ar +n}\fr{(\nu +n)}{1-q^{\nu+n}}=-1
\eqa
Equations (\r{rg}) and (\r{rp}) imply that 
\beq\lb{l}
\lim_{\nu \ar - n} \fr{\p (\nu )}{\G (\nu )} = (-1)^n q^{-n(n+1)/2} 
\fr{\log q}{1-q} \G (n+1).
\eeq

Before closing this section we like to present the $q\ar 1^-$ limit of 
$\p (\nu )$. We first rewrite it as 
\beq
\p (\nu)  =  \lim_{n\ar \infty } \sum_{l=1}^n 
(\log\fr{1-q^{l+1}}{1-q^l}+
\fr{q^{\nu + l-1}\log q} {1-q^{\nu+l-1}}) 
\eeq
Taking the $q\ar 1^-$ limit in the finite sum in the above formula
we have 
\beqa
\lim_{q\ar 1^-}\p (\nu)  =  \lim_{n\ar \infty } 
( \log (n+1) - \sum_{l=1}^n \fr{1}{\nu+l-1}) \nonumber \\
 =  \lim_{n\ar \infty } 
( \log (n+1) - \sum_{l=1}^{n+1} \fr{1}{l} + 
\sum_{l=1}^n (\fr{1}{l}-\fr{1}{\nu+l-1}) + \fr{1}{n+1}) \nonumber \\
= - C  + \sum_{l=0}^\infty  (\fr{1}{l+1}-\fr{1}{\nu+l}) = \psi (\nu ),
\eqa
where 
\beq
C = \lim_{n\ar \infty } ( \sum_{l=1}^{n} \fr{1}{l} -\log n) = -\psi (1)
\eeq
is the Euler number \cite{8}. It is then natural to define the 
q-Euler number 
\beq\lb{e}
C_q \equiv - \p (1).
\eeq

\vspace{1cm}
\noindent
{\large\bf 3.  q-Neumann  Function of Integer Order }

\vspace{2mm}
\noindent
The Hahn-Exton q-Bessel function of order $\nu$ is   defined as
\beq
J_\nu (x; q)= \sum_{k=0}^\infty \frac{(-1)^k q^{\fr{ k(k+1)}{2}} }
{\G (k+1)\G (k+\nu+1) } x^{2k+\nu}
\eeq
satisfies  the q-difference equation 
\beq\lb{d}
J_\nu (q^{1/2}x; q)+J_\nu (q^{-1/2}x; q)+ q^{-\nu /2}((1-q)^2 x^2 - q^\nu -1)
J_\nu (x; q)=0.
\eeq
For non-integer $\nu$ let us  define  the function
\beq\lb{n}
N_\nu (x; q)= \fr{ \cos (\pi\nu ) J_\nu (x; q)- q^{-\nu /2}J_{-\nu} 
(q^{-\nu /2}x; q)} { \sin (\pi\nu )}
\eeq
which is the second independent solution of the q-difference equation (\r{d}).
it  satisfies the same recurrence relations as the Hahn-Exton q-Bessel 
functions. For example we have  
\beq
q^{(\nu + 1)/2} N_{\nu +1} (q^{1/2}x; q) - N_{\nu +1} (x; q) 
=  (q-1)xN_\nu (x; q),
\eeq 
\beq
q^{\nu/2} N_\nu (q^{-1/2}x; q) - N_\nu (x; q) 
=  (q-1)xN_{\nu+1} (x; q),
\eeq 
Note that for integer $\nu =n$ one has the property 
\beq\lb{j}
J_{-n} (x; q) = (-1)^n q^{n/2}J_n (q^{n/2}x; q).
\eeq
Therefore to derive the form of (\r{n}) for integer order  
we use the L'Hospital rule 
\beq
\pi N_n (x; q)= \fr{d}{d\nu} J_\nu (x; q)\mid_{\nu =n} - (-1)^n 
\fr{d}{d\nu} ( q^{-\nu/2}J_{-\nu}(q^{-\nu/2}x; q) )\mid_{\nu =n}. 
\eeq
Making use of (\r{l}) and (\r{j}) for $n\in Z_+$ we get 
\beqa
\pi N_n (x; q)  =  2J_n (x; q)\log (q^{1/4}x) 
+\fr{\log q}{1-q} \sum_{k=0}^{n-1}\frac{\G (n-k) }
{\G (k+1) } x^{2k-n} \nonumber \\
 -  \sum_{k=0}^\infty \frac{(-1)^k q^{\fr{ k(k+1)}{2}} x^{2k+n} }
{\G (k+1)\G (k+n+1) } (\p (k+n+1)+\p (k+1) +k\log q). 
\eqa
Using the recurrence relation (\r{r}) and the definition (\r{e}) 
we arrive at 
\beqa\lb{n1}
\pi N_n (x; q)  =  2J_n (x; q)(\log (q^{1/4}x) +C_q)
+\fr{\log q}{1-q} \sum_{k=0}^{n-1}\frac{\G (n-k) }
{\G (k+1) } x^{2k-n} \nonumber \\
+\log q \sum_{k=1}^\infty \frac{(-1)^k q^{\fr{ k(k+1)}{2}} x^{2k+n} }
{\G (k+1)\G (k+n+1) } (\sum_{l=1}^{k+n}\fr{q^l}{1-q^l}+
\sum_{l=1}^k\fr{1}{1-q^l}) \nonumber \\ 
+\log q \fr{x^n}{\G (n+1)}\sum_{k=1}^n \fr{q^k}{1-q^k}.
\eqa
For $n=0$ we follow the similar steps and obtain
\beqa\lb{n2}
\pi N_0 (x; q) & = & 2J_0 (x; q)(\log (q^{1/4}x) +C_q) \nonumber \\
& + & \log q \sum_{k=1}^\infty \frac{(-1)^k q^{\fr{ k(k+1)}{2}} x^{2k} }
{\G (k+1)\G (k+1) } \sum_{l=1}^{k}\fr{1+q^l}{1-q^l}.
\eqa

It is obvious that (\r{n1}) and (\r{n2}) become the usual Neumann functions
in $q\ar 1^-$ limit. 
From the construction it is clear that the  q-Neumann functions 
of integer order satisfy the difference equation (\r{d}) and possess all 
recurrence relations satisfied by the Harh-Exton q-Bessel functions of 
integer order. We also have 
\beq
N_{-n} (x; q) = (-1)^n q^{n/2}N_n (q^{n/2}x; q).
\eeq

Before closing the section we note that several relations  
involving the q-Neumann  functions can be obtained from those 
of the Hahn-Exton q-Bessel functions. For example using the product 
formula \cite{10}
\beqa
 \sum_{s=-\infty}^\infty q^z J_x (q^{z/2}; q)J_{x-\nu} (q^{z/2}; q)
J_\nu (rq^{(y+\nu +z)/2}; q)= \nonumber \\
= J_0 (rq^{(x+y)/2}; q)J_\nu (rq^{(\nu+y)/2}; q)
\eqa
which is valid for $r$, $x$, $y$, $\nu \in C$;  $Re (x)>-1$, 
$\mid r\mid^2 q^{1+  Re(x)+ Re(y)} <1$ and $r\neq 0$ we obtain 
the product formula for  $x=-y=\nu/2$
\beq
\sum_{s=-\infty}^\infty q^z J_{\nu/2} (q^{z/2}; q)J_{-\nu/2} (q^{z/2}; q)
N_\nu (rq^{(\nu/4 +z/2)}; q)= J_0 (r; q)N_\nu (rq^{\nu/4}; q).
\eeq

\end{document}